\documentclass[12pt]{article}
\usepackage{amsfonts}
\usepackage{mathrsfs,dsfont,bbm}

\usepackage[latin1]{inputenc}
\usepackage{amsmath,amssymb}
\usepackage{latexsym}
\usepackage{epstopdf}
\usepackage{geometry}
\usepackage[active]{srcltx}
\usepackage[
bookmarks=true,         
bookmarksnumbered=true, 
colorlinks=true, pdfstartview=FitV, linkcolor=blue, citecolor=blue,
urlcolor=blue]{hyperref}

\newtheorem{theorem}{Theorem}[section]

\newtheorem{remark}{Remark}[section]

\numberwithin{equation}{section}
\begin{document}

\title{Existence and H\"{o}lder continuity conditions for self-intersection local time of Rosenblatt process}
\date{\today}

\author
{ Qian Yu, Guangjun Shen and Xiuwei Yin
\thanks{G. Shen is supported by National Natural Science Foundation of China (12071003).
X. Yin is supported by National Natural Science Foundation of China (11901005).} }

\maketitle

\begin{abstract}
\noindent We consider the existence and H\"{o}lder continuity conditions for the self-intersection local time of Rosenblatt process. Moreover, we  study the cases of intersection local time and collision local time, respectively.
\vskip.2cm \noindent {\bf Keywords:} Rosenblatt process; Self-intersection local time; Intersection local time; Collision local time; H\"{o}lder continuity.

\vskip.2cm \noindent {\it Subject Classification 2010: Primary 60G18; Secondary 60J55.}
\end{abstract}

\section{Introduction}\label{sec1}
Fractional Brownian motion (fBm) on $\mathbb{R}$  with Hurst parameter $H\in(0,1)$ is a  centered Gaussian process $B^H=\{B_t^H, ~t\geq0\}$ with covariance function given by
$$
\mathbb{E}\left[B^H_tB^H_s\right]=\frac{1}{2}\left[t^{2H}+s^{2H}-|t-s|^{2H}
\right].
$$
Note that $B_t^{\frac12}$ is a classical standard Brownian motion.
Let $D=\{(r,s): 0<r<s<t\}$. The self-intersection local time (SLT) of fBm was first investigated in Rosen \cite{Rosen 1987} and formally defined as
$$\alpha_t(y)=\int_{D}\delta(B^H_s-B^H_r-y)drds,$$
where $y\in\mathbb{R}$, $B^H$ is a fBm and $\delta$ is the Dirac delta function. Intuitively, $\alpha_t(0)$  measures the amount of time that the
process spends intersecting itself on the time interval $[0,t]$ and has been an important topic of the theory of stochastic process.
It was further investigated in Hu \cite{Hu 2001}, Hu and Nualart \cite{Hu 2005}.
In particular, Hu and Nualart \cite{Hu 2005} showed its existence whenever $H<1$. Moreover,
$\alpha_t(y)$ is H\"{o}lder continuous in time of any order strictly less than $1-H$ which can be derived from Xiao \cite{Xiao 1997}.

As we all know, fBm is a Gaussian process, and its characteristic function is in the form of exponential decay, which guarantees the finiteness of the moment  of SLT. However, for non-Gaussian processes, there is no such good characteristic function. This makes the corresponding results of local time not rich. To the best of our knowledge, the study of stable motion only appears in Yan \emph{et al.} \cite{Yan 2017}, Rosenblatt process case is only discussed in Bojdecki \emph{et al.} \cite{Boj2015} and Kerchev \emph{et al.} \cite{Ker2021}. In particular, Kerchev \emph{et al.} \cite{Ker2021} has studied the path properties of local time for Rosenblatt process, which gives us motivation to study SLT of Rosenblatt process in this paper. We formally define the SLT of Rosenblatt process at $y$ as
\begin{align*}
\widehat{\alpha}_t(y)=\int_{D}\delta(X^H_s-X^H_r-y)drds,
\end{align*}
where $X^H$ is a Rosenblatt process (the definition given in Section \ref{sec2}) with parameter $H\in(\frac12,1)$ and $\delta$ is the Dirac delta function.

Set
$$f_\varepsilon(x)=\frac1{\sqrt{2\pi\varepsilon}}e^{-\frac{|x|^2}{2\varepsilon}}=\frac1{2\pi}\int_{\mathbb{R}}e^{\iota px}e^{-\varepsilon \frac{|p|^2}{2}}dp.$$
Since the Dirac delta function $\delta$ can be approximated by $f_\varepsilon(x)$, we can approximate  $\widehat{\alpha}_t(y)$ by
\begin{equation}\label{sec1-eq1.2}
\widehat{\alpha}_{t,\varepsilon}(y)=\int_{D}f_\varepsilon(X^H_s-X^H_r-y)drds.
\end{equation}

If $\widehat{\alpha}_{t,\varepsilon}(y)$ converges to a random variable in $L^p$ as $\varepsilon\to0$, we denote the limit by $\widehat{\alpha}_t(y)$ and call it the SLT of $X^H$. Then  we can obtain the H\"{o}lder continuous both in time and space.
\begin{theorem}\label{sec1-DSLT}
If $\frac12<H<1$, then $\widehat{\alpha}_{t}(y)$ exists in $L^p$, for all $p\geq1$. Moreover,
$\widehat{\alpha}_{t}(y)$ is H\"{o}lder continuous in $y$ of any order strictly less than $\min\{1,\frac{1}{H}-1\}$
and H\"{o}lder continuous in $t$ of any order strictly  less than $1-H$, that is
\begin{equation*}
\left|\mathbb{E}\Big[\Big(\widehat{\alpha}_{t}(x)-\widehat{\alpha}_{t}(y)\Big)^n\Big]\right|\leq c\, |x-y|^{n\lambda},
\end{equation*}
where $\lambda<\min\{1,\frac{1}{H}-1\}$ and
\begin{equation*}
\left|\mathbb{E}\Big[\Big(\widehat{\alpha}_{t}(y)-\widehat{\alpha}_{\tilde{t}}(y)\Big)^n\Big]\right|\leq c\, |t-\tilde{t}|^{n\beta},
\end{equation*}
where $\beta<1-H$.
\end{theorem}

\begin{remark}\label{sec1-rem}
The condition $H>1/2$ is a necessary condition, which is introduced in the definition of Rosenblatt process and guarantees the finiteness of integral \eqref{sec3-DSLT-eq.4} in the proof of Theorem \ref{sec1-DSLT}.
Moreover, for $(d\geq2)$-dimensional Rosenblatt processes, we can define its SLT. However, the $L^p$ existence may not be found, because the finiteness of its $p$-th moment needs condition $Hd<1$, which is contrary to $H>1/2$.
\end{remark}

Similarly, we can define the intersection local time (ILT) of Rosenblatt process with parameter $H\in(1/2,1)$ and obtain its $L^p$ existence and H\"{o}lder  continuity. Define
\begin{align}\label{sec1-eq.DILT}
\widetilde{\alpha}_t(y)
=\int_{[0,t]^2}\delta^{(k)}(X^{H_1}_s-\widetilde{X}^{H_2}_r-y)drds,
\end{align}
where $X^{H_1}$ and $\widetilde{X}^{H_2}$ are two independent Rosenblatt processes. For fBm or more general Gaussian case, Guo \emph{et al.} \cite{Guo2017} and Hong and Xu \cite{HX2020,HX2021} study the existence and H\"{o}lder continuity. In this paper, we have the results of Rosenblatt case.

\begin{theorem}\label{sec1-DILT}
If $\frac{H_1H_2}{H_1+H_2}<1$, then $\widetilde{\alpha}_{t}(y)$ exists in $L^p$, for all $p\geq1$. Moreover, $\widetilde{\alpha}_{t}(y)$ satisfies the H\"{o}lder continuity:
\begin{equation*}
\left|\mathbb{E}\Big[\Big(\widetilde{\alpha}_{t}(x)-\widetilde{\alpha}_{t}(y)\Big)^n\Big]\right|\leq c\, |x-y|^{n\widetilde{\lambda}},
\end{equation*}
where $\widetilde{\lambda}<\min\{1,\frac{H_1+H_2}{2H_1H_2}-\frac{1}2\}$ and
\begin{equation*}
\left|\mathbb{E}\Big[\Big(\widetilde{\alpha}_{t}(y)-\widetilde{\alpha}_{\tilde{t}}(y)\Big)^n\Big]\right|\leq c\, |t-\tilde{t}|^{n\widetilde{\beta}},
\end{equation*}
where $\widetilde{\beta}<1-\frac{H_1H_2}{H_1+H_2}$.
\end{theorem}

Moreover, we can define the collision local time (CLT) of Rosenblatt process with parameter $H\in(1/2,1)$. Let
\begin{align}\label{sec1-eq.DCLT}
\overline{\alpha}_t(y)
&=\int_{[0,t]}\delta(X^{H_1}_s-\widetilde{X}^{H_2}_s-y)ds,
\end{align}
where $X^{H_1}$ and $\widetilde{X}^{H_2}$ and are two independent Rosenblatt processes. For Gaussian case, Jiang and Wang \cite{Jiang2009} and Shi \cite{Shi20} study the existence of CLT of bi-fBm. We now give the $L^p$ existence and H\"{o}lder  continuity of $\overline{\alpha}_{t}(y)$ as follows.

\begin{theorem}\label{sec1-DCLT}
If $\frac{H_1H_2}{H_1+H_2}<\frac12$, then $\overline{\alpha}_{t}(y)$ exists in $L^p$, for all $p\geq1$. Moreover, $\overline{\alpha}_{t}(y)$ satisfies the H\"{o}lder continuity:
\begin{equation*}
\left|\mathbb{E}\Big[\Big(\overline{\alpha}^{(k)}_{t}(x)-\overline{\alpha}^{(k)}_{t}(y)\Big)^n\Big]\right|\leq c\, |x-y|^{n\overline{\lambda}},
\end{equation*}
where $\overline{\lambda}<\min\{1,\frac{H_1+H_2}{4H_1H_2}-\frac{1}2\}$ and
\begin{equation*}
\left|\mathbb{E}\Big[\Big(\overline{\alpha}_{t}(y)-\overline{\alpha}_{\tilde{t}}(y)\Big)^n\Big]\right|\leq c\, |t-\tilde{t}|^{n\overline{\beta}},
\end{equation*}
where $\overline{\beta}<1-\frac{2H_1H_2}{H_1+H_2}$.
\end{theorem}

The topic of SLT of fBm has received a great deal of attention in recent decades. It was first studied by Rosen in \cite{Rosen 1987}, this has aroused the interest of many scholars in this research direction. The corresponding results are not only enriched from one-dimensional to multi-dimensional, but also extend the SLT itself to its derivative, including \cite{Das20+}, \cite{Hu 2001}, \cite{Hu 2005}, \cite{Jaramillo 2017}, \cite{Jaramillo 2019}, \cite{Jung 2014}, \cite{Jung 2015}, \cite{Rosen2005}, \cite{Yan 2015}, \cite{Yu20}--\cite{Yu21+} and references therein.

It should be noted that our results can not study the case of derivative of SLT for Rosenblatt process, because the Rosenblatt process can be rewrite as the second chaos (see \eqref{R-rep-2}), its characteristic function is similar to the chi square distribution, showing polynomial decay. If the derivative case is considered, we need to add the term of polynomial divergence in the integrand of \eqref{sec3-DSLT-eq.2}, which makes the finiteness of the integral cannot be guaranteed.

Although the H\"{o}lder continuous property of local time for Rosenblatt process has been proved in \cite{Ker2021} recently. As we all know, SLT and local time have different integral structures in form. In particular, for SLT, the second Rosenblatt process is dependent on the first Rosenblatt process in the integrand. So that technical methods given in \cite{Ker2021} can not be used directly here.

Moreover, we think the related results can be extended to the general cases. Since fBm ($q=1$) and Rosenblatt process ($q=2$) are both special Hermite process (with parameter $q\geq1$), they have many similar properties, such as the stationarity of increment, self similarity.
So, we believe that our results can be extended to Hermite processes, if Hermite processes are understood as $q$-th chaos and the corresponding spectral representation is obtained.

The paper has the following structure. Section \ref{sec2} contains some necessary preliminaries on Rosenblatt process. Section \ref{sec3} is to prove the main results. To be exact, we will split this section into three subsections to prove the three theorems given in Section \ref{sec1}.
Throughout this paper, if not mentioned otherwise, the letter $c$, with or without a subscript,
denotes a generic positive finite constant and may change from line to line.

\section{Preliminaries}\label{sec2}
In this section, we first give the definition of Rosenblatt process described in \cite{Tu08} and \cite{Tu13}, then give the result of spectral representation for Rosenblatt process obtained in \cite{Ker2021}.

The Hermite process (include Rosenblatt process) is an interesting class of self-similar processes with long range dependence, it is given as limits of the so called Non-Central Limit Theorem studied in Dobrushin and Major \cite{DM}, Taqqu \cite{Taqqu}. Let us briefly recall the general context.

Denote by  $H_j(x)$  the Hermite polynomial of order $j$ defined by $$H_j(x)=(-1)^je^{\frac{x^2}{2}}\frac{d^j}{dx^j}e^{-\frac{x^2}{2}},  \quad j=1,2,\cdots$$
with $H_0(x)=1$, and let the Borel function $g:\mathbb{R}\rightarrow\mathbb{R}$ satisfy $\mathbb{E}g(\xi_0)=0$, $\mathbb{E}g(\xi_0)^2<\infty$ and
$g(x)=\sum^{\infty}_{j=0}c_{j}H_{j}(x), ~c_{j}=\frac{1}{j!}E[g(\xi_0)H_{j}(\xi_0)]$.
The Hermite rank of $g$ is defined by $$q=\min\{j:c_j\neq0\}.$$ Clearly, $q\geq1$ since $\mathbb{E}g(\xi_0)=0$.

Let $g$ be a function of Hermite rank $q$ and let  $\{\xi_n, n\in\mathbb{N}\}$ be  stationary centered Gaussian sequence with $\mathbb{E}(\xi^2_{n})=1$ which exhibits long range dependence in the sense that the correlation function satisfies
$$
 r(n):=\mathbb{E}(\xi_{0}\xi_{n})=n^{\frac{2H-2}{q}}L(n),
$$
where $q\geq1$ is an integer, $H\in(\frac12,1)$ and $L$ is a slowly varying function at infinity. Then,
the Non Central Limit Theorem implies that the stochastic processes
$\frac{1}{n^H}\sum^{\lfloor nt\rfloor}_{j=1}g(\xi_j)$ ($\lfloor nt\rfloor$ denotes the integer part of $nt$)
converges, as $n\to\infty$, in the sense of finite dimensional distributions to the process
\begin{equation}\label{sec2-eq2.1}
X^{H,q}_t=c(H,q)\int_{\mathbb{R}^q}\int^t_0\left(\prod^q_{j=1}(s-y_j)_{+}^{-(\frac12+\frac{1-H}{q})}\right)dsB(dy_1)\cdot\cdot\cdot B(dy_q),
\end{equation}
where $x_+=\max\{x,0\}$ and the above integral is a Wiener-It\^{o} multiple integral with respect to the standard Brownian motion
$\{B(y), y\in\mathbb{R}\}$ excluding the diagonals $\{y_i=y_j\}, i\neq j$, $c(H,q)$ is a positive normalization constant depending
only on $H$ and $q$  such that $\mathbb{E}(X^{H,q}_1)^2=1$. The process $\{X^{H,q}_t, t\geq 0\}$ is called the \emph{Hermite process} of order $q$, it is
$H$ self-similar   and   has stationary increments. The class of Hermite processes includes fBm ($q=1$)
which is the only Gaussian process in this class. Since they are non-Gaussian ($q\geq 2$) and self-similar with stationary increments,
the Hermite processes can also be an input in models where self-similarity is observed in empirical data which appears to be non-Gaussian.
When $q=2$, the process  \eqref{sec2-eq2.1} is known as the \emph{Rosenblatt process}.

The Rosenblatt process $\{X^H_t, t\geq0\}$ admits the following stochastic representation (spectral representation):
\begin{equation}\label{R-rep-1}
X^H_t=\int_{\mathbb{R}^2}H_t(x,y)Z_G(dx)Z_G(dy),
\end{equation}
where the double Wiener-It\^{o} integral is taken over $x\neq y$, $Z_G(dx)$ is a complex-valued random white noise with control measure $G$ satisfying $G(dx)=|x|^{-H}dx$, and
$$H_t(x,y)=\frac{e^{\iota t(x+y)}-1}{\iota(x+y)}$$
is a complex valued Hilbert-Schmidt kernel with $H_t(x,y)=H_t(y,x)=\overline{H_t(-x,-y)}$ and
$$\int_{\mathbb{R}^2}|H_t(x,y)|^2G(dx)G(dy)<\infty.$$

In particular, by the spectral theorem, we can rewrite $X_t^H$ as an element in the second chaos
\begin{equation}\label{R-rep-2}
X^H_t\overset{Law}{=}\sum_{j=1}^\infty\lambda_j(Y_j^2-1),
\end{equation}
where $\{Y_j, j\geq1\}$ is a sequence of independent standard Gaussian random variable and $\{\lambda_j, j\geq1\}$ are the eigenvalues of the self-adjoint operator $A: L^2_G(\mathbb{R})\to L^2_G(\mathbb{R})$,
$$(Af)(x)=\int_{\mathbb{R}}H_t(x,-y)f(y)G(dy)=\int_{\mathbb{R}}H_t(x,-y)f(y)|x|^{-H}dy.$$
Moreover, we need condition $\sum_{j=1}^\infty\lambda_j^2<\infty$ to make sure \eqref{R-rep-2} converges.

\section{Proof of the main results}\label{sec3}

In this section, the proof of Theorems \ref{sec1-DSLT}, \ref{sec1-DILT} and \ref{sec1-DCLT} are taken into account. We will divide this
section into three parts and give the proof of the corresponding theorem in each part.

\subsection{The proof of Theorem \ref{sec1-DSLT}}\label{sec3.1}
\textbf{Existence in $L^p$.} We first prove that $\widehat{\alpha}_{t,\varepsilon}(y)$ converges in $L^p$ for $p\geq1$.
Since $|e^{-\iota p y}|=1$, we only need to consider the finiteness of $\left|\mathbb{E}(\widehat{\alpha}_{t,\varepsilon}(0))^p\right|$.

By \eqref{sec1-eq1.2}, $\widehat{\alpha}_{t,\varepsilon}(0)$ can be written as
\begin{align*}
\widehat{\alpha}_{t,\varepsilon}&:=\widehat{\alpha}_{t,\varepsilon}(0)=\int_0^t\int_0^sf_\varepsilon(X^H_s-X^H_r)drds\\
&=\frac{1}{2\pi}\int_{D}\int_{\mathbb{R}}e^{-\frac{\varepsilon}{2}p^2}e^{\iota p(X^H_s-X^H_r)}dpdrds,
\end{align*}
where $D=\{0<r<s<t\}$ and $\iota=\sqrt{-1}$.

For any integer $n\geq1$,
\begin{equation}\label{sec3-DSLT-eq.1}
\left|\mathbb{E}\left(\widehat{\alpha}_{t,\varepsilon}\right)^n\right|
\leq\frac{1}{(2\pi)^n}\int_{D^n}\int_{\mathbb{R}^{n}}\left|e^{-\frac{\varepsilon}{2}\sum_{i=1}^np_i^2}\right|\left|\mathbb{E}e^{\iota \sum_{i=1}^np_i(X^H_{s_i}-X^H_{r_i})}\right|dpdrds.
\end{equation}
We use the method of sample configuration as in Jung and Markowsky \cite{Jung 2015}.
Fix an ordering of the set $\{r_1, s_1, r_2, s_2, \cdots, r_n, s_n\}$, and let $l_1\leq l_2\leq \cdots \leq l_{2n}$ be a relabeling of the set $\{r_1, s_1, r_2, s_2, \cdots, r_n, s_n\}$. Let $u_1 \cdots u_{2n-1}$ be the proper linear combinations of the $p_{i}$'s so that
$$
\mathbb{E}\Big[e^{\iota \sum_{i=1}^np_i(X^H_{s_i}-X^H_{r_i})}\Big] = \mathbb{E}\Big[e^{\iota \sum_{i=1}^{2n-1}u_i(X^H_{\ell_{i+1}}-X^H_{\ell_i})}\Big].
$$

Therefore, we can rewrite \eqref{sec3-DSLT-eq.1} as
\begin{equation}\label{sec3-DSLT-eq.1+1}
\left|\mathbb{E}\left(\widehat{\alpha}_{t,\varepsilon}\right)^n\right|
\leq c\,\int_{E^n}\int_{\mathbb{R}^{n}}\left|\mathbb{E}e^{\iota \sum_{i=1}^{2n-1}u_i(X^H_{\ell_{i+1}}-X^H_{\ell_i})}\right|dpd\ell,
\end{equation}
where $E^n=\{0<\ell_1<\cdots<\ell_{2n}<t\}$.

By \eqref{R-rep-1}, we have
\begin{align*}
\sum_{i=1}^{2n-1}u_i(X^H_{\ell_{i+1}}-X^H_{\ell_i})=\int_{\mathbb{R}^2}\sum_{i=1}^{2n-1}u_i\frac{e^{\iota \ell_{i+1}(x+y)}-e^{\iota \ell_{i}(x+y)}}{\iota(x+y)}Z_G(dx)Z_G(dy),
\end{align*}
where the integral is taken over $x\neq y$. Define a operator $A_{\Delta\ell,u}$,
$$(A_{\Delta\ell,u}f)(x)=\int_{\mathbb{R}}\sum_{i=1}^{2n-1}u_i\frac{e^{\iota \ell_{i+1}(x+y)}-e^{\iota \ell_{i}(x+y)}}{\iota(x+y)}f(y)|y|^{-H}dy.$$
Thus, similar to \eqref{R-rep-2}, we have
\begin{equation}\label{R-sumrep-1}
\sum_{i=1}^{2n-1}u_i(X^H_{\ell_{i+1}}-X^H_{\ell_i})\overset{Law}{=}\sum_{j=1}^\infty\lambda_j(Y_j^2-1),
\end{equation}
where $\{Y_j, j\geq1\}$ is a sequence of independent standard Gaussian random variable and $\{\lambda_j, j\geq1\}$ are the eigenvalues of operator $A_{\Delta\ell,u}$.
$Y^2_j$  is a chi square distribution, and its characteristic function is $\mathbb{E}e^{\iota t Y_j^2}=(1-2\iota t)^{-\frac12}$. This gives
\begin{align*}
\left|\mathbb{E}e^{\iota \sum_{i=1}^np_i(X^H_{s_i}-X^H_{r_i})}\right|=\prod_{j=1}^\infty\frac{e^{-\iota\lambda_j}}{\sqrt{1-2\iota\lambda_j}}
=\prod_{j=1}^\infty(1+4\lambda_j^2)^{-\frac14}.
\end{align*}

Substitute the above equation into \eqref{sec3-DSLT-eq.1+1},
\begin{equation}\label{sec3-DSLT-eq.2}
\left|\mathbb{E}\left(\widehat{\alpha}_{t,\varepsilon}\right)^n\right|
\leq c\,\int_{E^n}\int_{\mathbb{R}^{n}}\prod_{j=1}^\infty(1+4\lambda_j^2)^{-\frac14}dpd\ell.
\end{equation}

Next, using the methods of Lemma 2.2 in \cite{Ker2021}, we obtain that
\begin{align}\label{sec3-DSLT-eq.3}
\lambda_{j}\geq c_{H}\max_{1\leq i\leq 2n-1}\{|u_i||\ell_{i+1}-\ell_{i}|^{H}\}\widetilde{\mu}_{j}^2
\end{align}
where $\widetilde{\mu}_n\sim \widetilde{c}_{H} n^{-\frac{H}{2}}$, and $c_{H}, ~\widetilde{c}_{H}>0$ are constants that only depends on $H$.

In fact, let $$B_{\Delta\ell,u}:=c_HK_{H/2}M_gK_{H/2}: ~L^2(\mathbb{R})\to L^2(\mathbb{R}),$$
where $g(x)=\sum_{i=1}^{2n-1}u_i\mathbf{1}_{[\ell_i,\ell_{i+1}]}(x)$, $M_g$ is the multiplication operator $(M_gf)(x)=g(x)f(x)$ and $K_{H/2}$ is a convolution operator defined via the Fourier transform $(\widehat{K_{H/2}f})(x)=|x|^{-H/2}\widehat{f}(x)$. Then we define operator
$$(Tf)(x)=|x|^{H/2}f(x).$$
Note that, $T: ~L^2(|y|^{-H}dy)\to L^2(\mathbb{R})$ is an  isometric isomorphism. Thus, the operator $A_{\Delta\ell,u}$ is isometrically isomorphic
to $V_{\Delta\ell,u}=TA_{\Delta\ell,u}T^{-1}: ~L^2(\mathbb{R})\to L^2(\mathbb{R})$ and satisfies
$$V_{\Delta\ell,u}f(x)=|x|^{H/2}\int_{\mathbb{R}}\sum_{i=1}^{2n-1}u_i\frac{e^{\iota \ell_{i+1}(x-y)}-e^{\iota \ell_{i}(x-y)}}{\iota(x-y)}f(y)|y|^{-H/2}dy.$$
The Fourier transform of $g$ $$(\mathcal{F}g)(x)=\left(\sum_{i=1}^{2n-1}u_i\mathcal{F}\mathbf{1}_{[\ell_i,\ell_{i+1}]}\right)(x)=\sum_{i=1}^{2n-1}u_i\frac{e^{-\iota\ell_{i+1}x}-e^{-\iota\ell_{i}x}}{-\iota x}$$
gives
\begin{align*}
\mathcal{F}^3(K_{H/2}M_gK_{H/2}f)(x)&=2\pi\mathcal{F}(K_{H/2}M_gK_{H/2}f)(-x)\\
&=2\pi|x|^{-H/2}\int_{\mathbb{R}}\sum_{i=1}^{2n-1}u_i\frac{e^{-\iota\ell_{i+1}x}-e^{-\iota\ell_{i}x}}{-\iota x}\widehat{f}(-y)|y|^{-H/2}dy\\
&=|x|^{-H/2}\int_{\mathbb{R}}\sum_{i=1}^{2n-1}u_i\frac{e^{-\iota\ell_{i+1}x}-e^{-\iota\ell_{i}x}}{-\iota x}(\mathcal{F}^3f)(y)|y|^{-H/2}dy\\
&=(V_{\Delta\ell,u}\mathcal{F}^3f)(x).
\end{align*}
Hence, $A_{\Delta\ell,u}$ and $c_HK_{H/2}M_gK_{H/2}=:B_{\Delta\ell,u}$ are unitarily equivalent and have the same eigenvalues.

By the Lemma 5.4 in \cite{Ker2021},  we can find
$$\mu_j(M_{\mathbf{1}_{[\ell_i,\ell_{i+1}]}}K_{H/2}M_{\mathbf{1}_{[\ell_i,\ell_{i+1}]}})\geq c_{H}|\ell_{i+1}-\ell_i|^{H/2}\widetilde{\mu}_j$$
and the eigenvalues $\mu_j$ of $B_{\Delta\ell,u}$ greater than
\begin{align*}
\mu_j(K_{H/2}M_gK_{H/2})
&\geq\mu_j\left[(M_{g}K_{H/2}M_{g})^2\right]\\
&=\left(\mu_j(M_{g}K_{H/2}M_{g})\right)^2\nonumber\\
&\geq c_{H}\max_{1\leq i\leq2n-1}|u_i||\ell_{i+1}-\ell_i|^{H}\widetilde{\mu}_j^2,
\end{align*}
where $\widetilde{\mu}_j\sim j^{-H/2}$.

Substituting \eqref{sec3-DSLT-eq.3} into \eqref{sec3-DSLT-eq.2},
\begin{align}\label{sec3-DSLT-eq.4}
\left|\mathbb{E}\left(\widehat{\alpha}_{t,\varepsilon}\right)^n\right|
&\leq c\,\int_{E^n}\int_{\mathbb{R}^n}\prod_{j=1}^\infty\left(1+c_H\max_{1\leq i\leq2n-1}|u_i|^2|\ell_{i+1}-\ell_i|^{2H}\widetilde{\mu}_j^4\right)^{-\frac14}dpd\ell\nonumber\\
&\leq c\,\int_{E^n}\int_{\mathbb{R}^n}\prod_{j=1}^\infty\left(1+c_H\max_{i\in\mathcal{A}}|u_i|^2|\ell_{i+1}-\ell_i|^{2H}\widetilde{\mu}_j^4\right)^{-\frac14}dpd\ell\nonumber\\
&\leq c\,\int_{E^n}\int_{\mathbb{R}^n}\exp(-c_H\max_{i\in\mathcal{A}}|u_i||\ell_{i+1}-\ell_i|^{H})dpd\ell\nonumber\\
&\leq c_H\,|J|\int_{E^n}\int_{\mathbb{R}^n}
\left[\prod_{i\in\mathcal{A}}\exp(-c_H|u_i|^{\frac1H}|\ell_{i+1}-\ell_i|)\right]^{\frac1n}dud\ell,
\end{align}
where $\mathcal{A}$ is a subset of $\{1,2,\cdots, 2n-1\}$ such that the set $\{u_i, i\in\mathcal{A}\}$ spans
$\{p_1, p_2,\cdots, p_n\}$, $|J|$ is the Jacobian determinant of changing variables $(p_1, p_2,\cdots, p_n)$ to $\{u_i, i\in\mathcal{A}\}$, we use the inequality (see in page 518 line -5 of \cite{Ker2021})
\begin{align}\label{sec3-DSLT-eq.ref}
\prod_{j=1}^\infty\left(1+c_H\max_{i}|u_i|^2|\ell_{i+1}-\ell_i|^{2H}\widetilde{\mu}_j^4\right)^{-\frac14}\leq c\,\exp(-c_H\max_{i}|u_i|^{\frac1H}|\ell_{i+1}-\ell_i|)
\end{align}
with $H\in(\frac12, 1)$, in the third inequality.

It is easy to see
$$\int_{\mathbb{R}}\exp(-c_H|u|^{\frac1H}|\Delta\ell_{i+1}|)du\leq c_H|\Delta\ell_{i+1}|^{-H}.$$
Thus, the integral with respect to $du$ in \eqref{sec3-DSLT-eq.4} is
\begin{align*}
\int_{\mathbb{R}^n}
\left[\prod_{i\in\mathcal{A}}\exp(-c_H|u_i|^{\frac1H}|\ell_{i+1}-\ell_i|)\right]^{\frac1n}du\leq c_{H,n}\prod_{i\in\mathcal{A}}|\ell_{i+1}-\ell_i|^{-H}.
\end{align*}
This gives
\begin{align}\label{sec3-DSLT-eq.en}
\left|\mathbb{E}\left(\widehat{\alpha}_{t,\varepsilon}\right)^n\right|
\leq c_{H,n}\int_{E^n}\prod_{i\in\mathcal{A}}|\ell_{i+1}-\ell_i|^{-H}d\ell,
\end{align}
which is finite since $H<1$.

Next, we need to prove that $\left\{\widehat{\alpha}_{t,\varepsilon}, \varepsilon>0\right\}$ is a Cauchy sequence. For any $\varepsilon, \eta>0$,
\begin{align*}
\left|\mathbb{E}\left(\widehat{\alpha}_{t,\varepsilon}-\widehat{\alpha}_{t,\eta}\right)^n\right|
&\leq\frac{1}{(2\pi)^n}\int_{E^{n}}\int_{\mathbb{R}^n}\prod_{i=1}^n\left|e^{-\frac{\varepsilon}2|p_i|^2}-e^{-\frac{\eta}2|p_i|^2}\right|\prod_{j=1}^\infty(1+4\lambda_{j}^2)^{-\frac14}dpd\ell.
\end{align*}
By the dominated convergence theorem and
\begin{align*}
\int_{E^{n}}\int_{\mathbb{R}^n}
\prod_{j=1}^\infty(1+4\lambda_{j}^2)^{-\frac14}dpd\ell<\infty.
\end{align*}
So, we can see that $\widehat{\alpha}_{t,\varepsilon}$ exists in $L^p$, for all $p\geq1$, under the condition $1/2<H<1$.

\textbf{H\"{o}lder continuity of space variable.}
For the H\"{o}lder continuity of $\widehat{\alpha}_{t}(y)$ in variable $y$, define $x, y\in\mathbb{R}$,
\begin{align*}
&\left|\mathbb{E}(\widehat{\alpha}_{t,\varepsilon}(x)-\widehat{\alpha}_{t,\varepsilon}(y))^n\right|\\
&\leq\frac1{(2\pi)^{n}}\int_{E^{n}}\int_{\mathbb{R}^{n}}\prod_{j=1}^n\left(e^{-\frac{\varepsilon}{2}|p_j|^2}\left|e^{-\iota p_j x}-e^{-\iota p_jy}\right|\right)\prod_{j=1}^\infty(1+4\lambda_{j}^2)^{-\frac14}dpd\ell.
\end{align*}
Note that $$\left|e^{-\iota p_jx}-e^{-\iota p_jy}\right|\leq c_{\lambda}|x-y|^\lambda|p_j|^\lambda, ~~\lambda\in[0,1].$$
Fix $j$, and let $j_1$ to be the smallest value such that $u_{j_1}$ contains $p_{j}$ as a term and then choose $j_2$ to be
the smallest value strictly larger than $j_1$ such that $u_{j_2}$ does not contain $p_{j}$ as a term. Then $p_j = u_{j_1} - u_{j_1-1} = u_{j_2-1} - u_{j_2}$. We can see that, with the convention that $u_{0}=u_{2n}=0$,
\begin{align*}
|p_j|^\lambda=|u_{j_1}-u_{j_1-1}|^{\frac {\lambda}2}|u_{j_2-1}-u_{j_2}|^{\frac {\lambda}2} \leq c\left(|u_{j_1}|^{\frac {\lambda}2}+|u_{j_1-1}|^{\frac {\lambda}2}\right)\left(|u_{j_2-1}|^{\frac {\lambda}2}+|u_{j_2}|^{\frac {\lambda}2} \right).
\end{align*}
Thus,
\begin{align*}
\prod_{i=1}^n|p_i|^\lambda &\leq c \prod_{i=1}^{2n} (|u_{i}|^{\frac {\lambda}2}+|u_{i-1}|^{\frac {\lambda}2})\\
&=c\,\sum_{S_1}\prod_{i=1}^{2n}(|u_{i}|^{\frac {\lambda}2\gamma_{i}}|u_{i-1}|^{\frac {\lambda}2\overline{\gamma_{i}}})\\
&\leq c\,\sum_{S_2}\prod_{i=1}^{2n-1}(|u_{i}|^{\frac {\lambda}2\alpha_{i}}),
\end{align*}
where $$S_1=\left\{\gamma_{i}, ~\overline{\gamma_{i}}: ~\gamma_{i}\in\{0,1\}, ~\gamma_{i}+\overline{\gamma_{i}}=1,  ~i=1, \cdots, 2n\right\}$$
and
$$
S_2=\left\{\alpha_{i}: ~\alpha_{i}\in\{0, 1, 2\},  ~i=1, \cdots 2n-1\right\}.
$$
Then
\begin{align*}
\left|\mathbb{E}(\widehat{\alpha}_{t,\varepsilon}(x)-\widehat{\alpha}_{t,\varepsilon}(y))^n\right|\leq c\,|x-y|^\lambda\sum_{S_2}\int_{E^{n}}\int_{\mathbb{R}^{n}}\prod_{i=1}^{2n-1}|u_{i}|^{\frac{\lambda}{2}\alpha_i}
\exp\prod_{j=1}^\infty(1+4\lambda_{j}^2)^{-\frac14}dpd\ell.
\end{align*}

Similar to \eqref{sec3-DSLT-eq.4},
\begin{align*}
&\int_{\mathbb{R}^n}\prod_{i=1}^{2n-1}|u_{i}|^{\frac{\lambda}{2}\alpha_i}\prod_{j_1=1}^\infty(1+4\lambda_{j_1}^2)^{-\frac14}dp\\
&\leq c\,\int_{\mathbb{R}^n}\prod_{i=1}^{2n-1}|u_{i}|^{\frac{\lambda}{2}\alpha_i}\prod_{j=1}^\infty\left(1+c_H\max_{1\leq i\leq2n-1}|u_i|^2|\ell_{i+1}-\ell_i|^{2H}\widetilde{\mu}_j^4\right)^{-\frac14}dpd\ell\\
&\leq c\,\int_{\mathbb{R}^n}\prod_{i=1}^{2n-1}|u_{i}|^{\frac{\lambda}{2}\alpha_i}\exp(-c_H\max_{1\leq i\leq2n-1}|u_i|^{\frac1H}|\ell_{i+1}-\ell_i|)dp\\
&\leq c\,\int_{\mathbb{R}^n}\prod_{i\in\mathcal{A}^C}\left(|u_{i}|^{\frac{\lambda}{2}\alpha_i}\exp(-\frac{c_H}{2n}|u_i|^{\frac1H}|\ell_{i+1}-\ell_i|)\right)\\
&\qquad\times \prod_{i\in\mathcal{A}}\left(|u_{i}|^{\frac{\lambda}{2}\alpha_i}\exp(-\frac{c_H}{2n}|u_i|^{\frac1H}|\ell_{i+1}-\ell_i|)\right)dp\\
&\leq c\,\prod_{i\in\mathcal{A}^C}|\ell_{i+1}-\ell_i|^{-\frac{\lambda\alpha_i}{2}H}|J|\int_{\mathbb{R}^n}\prod_{i\in\mathcal{A}}\left(|u_{i}|^{\frac{\lambda}{2}\alpha_i}\exp(-\frac{c_H}{2n}|u_i|^{\frac1H}|\ell_{i+1}-\ell_i|)\right)du\\
&\leq c_{H,n}\prod_{i\in\mathcal{A}^C}|\ell_{i+1}-\ell_i|^{-\frac{\lambda\alpha_i}{2}H}\prod_{i\in\mathcal{A}}|\ell_{i+1}-\ell_i|^{-\frac{\lambda\alpha_i}{2}H-H},
\end{align*}
where $\mathcal{A}^C$ denotes the complement of $\mathcal{A}$ in $\{1,2,\cdots,2n-1\}$, $|J|$ is the Jacobian determinant of changing variables $(p_1, p_2,\cdots, p_n)$ to $\{u_i, i\in\mathcal{A}\}$ and we use the inequality
$$|u_{i}|^{\frac{\lambda}{2}\alpha_i}\exp(-\frac{c_H}{2n}|u_i|^{\frac1H}|\ell_{i+1}-\ell_i|)\leq c_{H,n}\,|\ell_{i+1}-\ell_i|^{-\frac{\lambda\alpha_i}{2}H}$$
in the last second inequality.

Thus,
\begin{align*}
\left|\mathbb{E}(\widehat{\alpha}_{t,\varepsilon}(x)-\widehat{\alpha}_{t,\varepsilon}(y))^n\right|\leq c_{H,n}\,|x-y|^\lambda\sum_{S_2}\int_{E^{n}}\prod_{i\in\mathcal{A}^C}|\ell_{i+1}-\ell_i|^{-\frac{\lambda\alpha_i}{2}H}\prod_{i\in\mathcal{A}}|\ell_{i+1}-\ell_i|^{-\frac{\lambda\alpha_i}{2}H-H}d\ell.
\end{align*}
which is finite with $H(\lambda+1)<1$, since $\alpha_{i}\in\{0, 1, 2\}$.

So,
we need the condition $H(\lambda+1)<1$ to make sure
\begin{align*}
&\left|\mathbb{E}(\widehat{\alpha}_{t}(x)-\widehat{\alpha}_{t}(y))^n\right|\\
&\leq c_n\left|\mathbb{E}(\widehat{\alpha}_{t}(x)-\widehat{\alpha}_{t,\varepsilon}(x))^n\right|
+c_n\left|\mathbb{E}(\widehat{\alpha}_{t,\varepsilon}(x)-\widehat{\alpha}_{t,\varepsilon}(y))^n\right|\\
&\qquad+c_n\left|\mathbb{E}(\widehat{\alpha}_{t,\varepsilon}(y)-\widehat{\alpha}_{t}(y))^n\right|\\
&\leq c_{n,H,t}\,|x-y|^{n\lambda}.
\end{align*}

Hence, from the Kolmogorov continuity criterion, the H\"{o}lder continuity of $\widehat{\alpha}_{t}(y)$ in space variable of any order $\lambda$ strictly less than $\min\{1,\frac1H-1\}$.

\textbf{H\"{o}lder continuity of time variable.}
For the H\"{o}lder continuity of $\widehat{\alpha}_{t}(y)$ in time variables $t$. Without loss of generality, we assume that $t<\widehat{t}$ and let $\widehat{D}=\{(r,s): ~0<r<s<\widehat{t}\}$.

Then
\begin{align*}
&\left|\mathbb{E}(\widehat{\alpha}_{t,\varepsilon}(y)-\widehat{\alpha}_{\widehat{t},\varepsilon}(y))^n\right|\\
&\leq\frac{1}{(2\pi)^n}\int_{(\widehat{D}\setminus D)^n}\int_{\mathbb{R}^{n}}\left|\mathbb{E}e^{\iota \sum_{i=1}^np_i(X^H_{s_i}-X^H_{r_i})}\right|dpdrds\\
&\leq c\,\int_{[t,\widehat{t}]^n}\int_{[0,s_1]\times\cdots\times[0,s_n]}\int_{\mathbb{R}^{n}}\left|\mathbb{E}e^{\iota \sum_{i=1}^np_i(X^H_{s_i}-X^H_{r_i})}\right|dpdrds\\
&\leq c\,\int_{\widehat{D}^n}\prod_{i=1}^n\mathbf{1}_{[t,\widehat{t}]}(s_i)\int_{\mathbb{R}^{n}}\left|\mathbb{E}e^{\iota \sum_{i=1}^np_i(X^H_{s_i}-X^H_{r_i})}\right|dpdrds\\
&\leq c\,|\widehat{t}-t|^{n\beta}\left(\int_{\widehat{D}^n}\left(\int_{\mathbb{R}^{n}}\left|\mathbb{E}e^{\iota \sum_{i=1}^np_i(X^H_{s_i}-X^H_{r_i})}\right|dp\right)^{\frac1{1-\beta}}drds\right)^{1-\beta}\\
&=:c\,|\widehat{t}-t|^{n\beta}\Lambda,
\end{align*}
where we use the H\"{o}lder's inequality in the last inequality with $\beta<1-H$.

Using the similar methods as in \eqref{sec3-DSLT-eq.en}, $\Lambda$ is bounded by
\begin{align*}
\left|\mathbb{E}\left(\widehat{\alpha}_{t,\varepsilon}\right)^n\right|
\leq c_{H,n}\left(\int_{E^n}\prod_{i\in\mathcal{A}}|\ell_{i+1}-\ell_i|^{-H/(1-\beta)}d\ell\right)^{1-\beta},
\end{align*}
since $1-H/(1-\beta)>0$.

This gives
\begin{align*}
&\left|\mathbb{E}(\widehat{\alpha}_{t}(x)-\widehat{\alpha}_{\widetilde{t}}(x))^n\right|\\
&\leq c_n\left|\mathbb{E}(\widehat{\alpha}_{t}(x)-\widehat{\alpha}_{t,\varepsilon}(x))^n\right|
+c_n\left|\mathbb{E}(\widehat{\alpha}_{t,\varepsilon}(x)-\widehat{\alpha}_{\widetilde{t},\varepsilon}(x))^n\right|\\
&\qquad+c_n\left|\mathbb{E}(\widehat{\alpha}_{\widetilde{t},\varepsilon}(x)-\widehat{\alpha}_{\widetilde{t}}(x))^n\right|\\
&\leq c_{n,H_1,H_2,t}\,|\widehat{t}-t|^{n\beta}
\end{align*}
with $\beta<1-H$.  This completes the proof.

\subsection{The proof of Theorem \ref{sec1-DILT}}\label{sec3.2}
In this subsection, we will consider the case of ILT of Rosenblatt process.

By the definition \eqref{sec1-eq.DILT}, we let
\begin{align*}
\widetilde{\alpha}_{t,\varepsilon}:&=\widetilde{\alpha}_{t,\varepsilon}(0)=\int_0^t\int_0^tf_\varepsilon(X^{H_1}_s-\widetilde{X}^{H_2}_r)drds\\
&=\frac{1}{2\pi}\int_{[0,t]^2}\int_{\mathbb{R}}e^{-\frac{\varepsilon}{2}p^2}e^{\iota p(X^{H_1}_s-\widetilde{X}^{H_2}_r)}dpdrds,
\end{align*}
where $X^{H_1}$ and $\widetilde{X}^{H_2}$ are two independent Rosenblatt processes.

\textbf{Existence in $L^p$.}
For positive integer $n\geq1$,
\begin{align}\label{sec3-DILT-eq.1}
\left|\mathbb{E}\left(\widetilde{\alpha}_{t,\varepsilon}\right)^n\right|
\leq\frac{1}{(2\pi)^n}\int_{[0,t]^{2n}}\int_{\mathbb{R}^n}\left|\mathbb{E}e^{\iota \sum_{i=1}^np_i(X^{H_1}_s-\widetilde{X}^{H_2}_r)}\right|dpdrds.
\end{align}

Define two operators $A_{s,p}$ and $\widetilde{A}_{r,p}$,
$$(A_{s,p}f)(x)=\int_{\mathbb{R}}\sum_{i=1}^np_i\frac{e^{\iota s_i(x+y)}-1}{\iota(x+y)}f(y)|y|^{-H_1}dy$$
and
$$(\widetilde{A}_{r,p}f)(x)=\int_{\mathbb{R}}\sum_{i=1}^np_i\frac{e^{\iota r_i(x+y)}-1}{\iota(x+y)}f(y)|y|^{-H_2}dy.$$
Similar to \eqref{R-sumrep-1}, we have
\begin{equation*}
\sum_{i=1}^np_i(X^{H_1}_{s_i}-\widetilde{X}^{H_2}_{r_i})\overset{Law}{=}\sum_{j=1}^\infty(\lambda_j-\widetilde{\lambda}_j)(Y_j^2-1),
\end{equation*}
where $\{Y_j, j\geq1\}$ is a sequence of independent standard Gaussian random variable, $\{\lambda_j, j\geq1\}$ and $\{\widetilde{\lambda}_j, j\geq1\}$ are the eigenvalues of operator $A_{s,p}$ and $\widetilde{A}_{r,p}$, respectively. Then
\begin{align*}
\left|\mathbb{E}e^{\iota \sum_{i=1}^np_i(X^{H_1}_{s_i}-\widetilde{X}^{H_2}_{r_i})}\right|
&=\prod_{j_1=1}^\infty\frac{e^{-\iota\lambda_{j_1}}}{\sqrt{1-2\iota\lambda_{j_1}}}\prod_{j_2=1}^\infty\frac{e^{\iota\widetilde{\lambda}_{j_2}}}{\sqrt{1+2\iota\widetilde{\lambda}_{j_2}}}\\
&=\prod_{j_1=1}^\infty(1+4\lambda_{j_1}^2)^{-\frac14}\prod_{j_2=1}^\infty(1+4\widetilde{\lambda}_{j_2}^2)^{-\frac14}.
\end{align*}

Substituting the above equation into \eqref{sec3-DILT-eq.1},
\begin{align}\label{sec3-DILT-eq.2}
\left|\mathbb{E}\left(\widetilde{\alpha}_{t,\varepsilon}\right)^n\right|
&\leq\frac{1}{(2\pi)^n}\int_{[0,t]^{2n}}\int_{\mathbb{R}^n}\prod_{j_1=1}^\infty(1+4\lambda_{j_1}^2)^{-\frac14}\prod_{j_2=1}^\infty(1+4\widetilde{\lambda}_{j_2}^2)^{-\frac14}dpdrds\nonumber\\
&\leq c\,\int_{[0,t]^{2n}}\left(\int_{\mathbb{R}^n}\prod_{j_1=1}^\infty(1+4\lambda_{j_1}^2)^{-\frac14\frac{H_1+H_2}{H_2}}dp\right)^{\frac{H_2}{H_1+H_2}}\nonumber\\
&\qquad\times\left(\int_{\mathbb{R}^n}\prod_{j_2=1}^\infty(1+4\widetilde{\lambda}_{j_2}^2)^{-\frac14\frac{H_1+H_2}{H_1}}dp\right)^{\frac{H_1}{H_1+H_2}}drds\nonumber\\
&\leq c_n\, \int_{\widetilde{D}^{2n}}\left(\int_{\mathbb{R}^n}\prod_{j_1=1}^\infty(1+4\lambda_{j_1}^2)^{-\frac14}dp\right)^{\frac{H_2}{H_1+H_2}}\nonumber\\
&\qquad\times\left(\int_{\mathbb{R}^n}\prod_{j_2=1}^\infty(1+4\widetilde{\lambda}_{j_2}^2)^{-\frac14}dp\right)^{\frac{H_1}{H_1+H_2}}drds,
\end{align}
where $\widetilde{D}=\{0<s_1<s_2<\cdots<s_n<t\}$.

Similar to the case of SLT, note that $A_{s,u}$ and $c_{H_1}K_{H_1/2}M_{\widetilde{g}}K_{H_1/2}$ are unitarily equivalent and have the same eigenvalues,
where
$$\widetilde{g}(x)=\sum_{i=1}^np_i\mathbf{1}_{[0,s_i]}(x)=\sum_{i=1}^n(p_i-p_{i-1})\mathbf{1}_{[s_{i-1},s_i]}(x)$$
and $0=s_0<s_1<\cdots<s_n$ with the convention $p_0=0$.

By \eqref{sec3-DSLT-eq.3}, the $j_1$-th singular value $\lambda_{j_1}$ of $A_{s,p}$ satisfies
\begin{align*}
\lambda_{j_1}=\mu_{j_1}(K_{H_1/2}M_{\widetilde{g}}K_{H_1/2})\geq c_{H_1}\max_{1\leq j\leq n}|p_j-p_{j-1}||s_j-s_{j-1}|^{H_1}\}\widetilde{\mu}_{j_1}^2,
\end{align*}
where $\widetilde{\mu}_n\sim \widetilde{c}_{H_1} n^{-\frac{H_1}{2}}$, and $c_{H_1}, ~\widetilde{c}_{H_1}>0$ are constants that only depends on $H_1$.

This gives
\begin{align}\label{sec3-DILT-eq.3}
&\int_{\mathbb{R}^n}\prod_{j_1=1}^\infty(1+4\lambda_{j_1}^2)^{-\frac14}dp\nonumber\\
&\leq c\,\int_{\mathbb{R}^n}\prod_{j_1=1}^\infty\left(1+c_{H_1}\,(\max_{1\leq j\leq n}|p_j-p_{j-1}||s_j-s_{j-1}|^{H_1})^2\widetilde{\mu}_{j_1}^4\right)^{-\frac14}dp\nonumber\\
&\leq c_{n}\,\int_{\mathbb{R}^n}\prod_{j_1=1}^\infty\left(1+c_{H_1}(\max_{1\leq j\leq n}|s_j-s_{j-1}|^{H_1}|\xi_j|)^2\widetilde{\mu}_{j_1}^4\right)^{-\frac14}d\xi,
\end{align}
where we make the change of variables $\xi_j=\sum_{i=j}^np_i$ for $j=1,2,\cdots, n$ with the convention $\xi_{n+1}=0$ and $\max_{j}|2\xi_j-\xi_{j-1}-\xi_{j+1}|\geq c_n\,\max_j|\xi_j|$ in the last inequality.

For the integrand function in \eqref{sec3-DILT-eq.3}, similar to \eqref{sec3-DSLT-eq.ref}, we have
\begin{align}\label{sec3-DILT-eq.ref}
\prod_{j_1=1}^\infty\left(1+c_{H_1}(\max_{1\leq j\leq n}|s_j-s_{j-1}|^{H_1}|\xi_j|)^2\widetilde{\mu}_{j_1}^4\right)^{-\frac14}\leq c\,\exp(-c_{H_1}\max_{1\leq j\leq n}|s_j-s_{j-1}||\xi_j|^{\frac1{H_1}}).
\end{align}

Substituting \eqref{sec3-DILT-eq.ref} into \eqref{sec3-DILT-eq.3},
\begin{align}\label{sec3-DILT-eq.5}
\int_{\mathbb{R}^n}\prod_{j_1=1}^\infty(1+4\lambda_{j_1}^2)^{-\frac14}dp
&\leq c\,\int_{\mathbb{R}^n}e^{-c_{H_1}\max_{1\leq j\leq n}|s_j-s_{j-1}||\xi_j|^{\frac1{H_1}}}d\xi\nonumber\\
&\leq c\,\int_{\mathbb{R}^n}\left(e^{-c_{H_1}\sum_{j=1}^n|s_j-s_{j-1}||\xi_j|^{\frac1{H_1}}}\right)^{\frac1n}d\xi\nonumber\\
&\leq c_{H_1,n}\prod_{j=1}^n|s_j-s_{j-1}|^{-H_1}.
\end{align}

Similarly, we have the same inequality with respect to $\widetilde{\lambda}_{j_2}$,
\begin{align}\label{sec3-DILT-eq.6}
\int_{\mathbb{R}^n}\prod_{j_2=1}^\infty(1+4\widetilde{\lambda}_{j_2}^2)^{-\frac14}dp
\leq c_{H_2,n}\prod_{j=1}^n|s_j-s_{j-1}|^{-H_2}.
\end{align}

Together \eqref{sec3-DILT-eq.2}, \eqref{sec3-DILT-eq.5} and \eqref{sec3-DILT-eq.6}, we have
\begin{align}\label{sec3-DILT-eq.7}
\left|\mathbb{E}\left(\widetilde{\alpha}_{t,\varepsilon}^{(k)}\right)^n\right|
&\leq c_n\left(\int_{[0,t]^n}\prod_{j=1}^n|r_j-r_{j-1}|^{-\frac{H_1H_2}{H_1+H_2}}dr\right)^2,
\end{align}
which is finite under condition $\frac{H_1H_2}{H_1+H_2}<1$.

Next, we need to prove that $\left\{\widetilde{\alpha}_{t,\varepsilon}, \varepsilon>0\right\}$ is a Cauchy sequence. For any $\varepsilon, \eta>0$,
\begin{align*}
\left|\mathbb{E}\left(\widetilde{\alpha}_{t,\varepsilon}-\widetilde{\alpha}_{t,\eta}\right)^n\right|
&\leq\frac{1}{2\pi}\int_{[0,t]^{2n}}\int_{\mathbb{R}^n}\prod_{i=1}^n\left|e^{-\frac{\varepsilon}2|p_i|^2}-e^{-\frac{\eta}2|p_i|^2}\right|\\
&\qquad\quad\times\prod_{j_1=1}^\infty(1+4\lambda_{j_1}^2)^{-\frac14}\prod_{j_2=1}^\infty(1+4\widetilde{\lambda}_{j_2}^2)^{-\frac14}dpdrds.
\end{align*}
By the dominated convergence theorem and
\begin{align*}
\int_{[0,t]^{2n}}\int_{\mathbb{R}^n}
\prod_{j_1=1}^\infty(1+4\lambda_{j_1}^2)^{-\frac14}\prod_{j_2=1}^\infty(1+4\widetilde{\lambda}_{j_2}^2)^{-\frac14}dpdrds<\infty
\end{align*}
with $\frac{H_1H_2}{H_1+H_2}<1$. So, we can see that $\widetilde{\alpha}_{t,\varepsilon}$ exists in $L^p$, for all
$p\geq1$, under the condition $H_1H_2<H_1+H_2$.

\textbf{H\"{o}lder continuity of space variable.}
For $x,y\in\mathbb{R}$,
\begin{align*}
&\left|\mathbb{E}(\widetilde{\alpha}_{t,\varepsilon}(x)-\widetilde{\alpha}_{t,\varepsilon}(y))^n\right|\\
&\leq\frac1{(2\pi)^{n}}\int_{[0,t]^{2n}}\int_{\mathbb{R}^{n}}\prod_{j=1}^n\left(e^{-\frac{\varepsilon}{2}|p_j|^2}\left|e^{-\iota p_j x}-e^{-\iota p_jy}\right|\right)\\
&\qquad\times\prod_{j_1=1}^\infty(1+4\lambda_{j_1}^2)^{-\frac14}\prod_{j_2=1}^\infty(1+4\widetilde{\lambda}_{j_2}^2)^{-\frac14}dpdrds.
\end{align*}
Since $$\left|e^{-\iota p_jx}-e^{-\iota p_jy}\right|\leq c_{\lambda}|x-y|^\lambda|p_j|^\lambda, ~~\lambda\in[0,1],$$
Then
\begin{align*}
&\left|\mathbb{E}(\widetilde{\alpha}_{t,\varepsilon}(x)-\widetilde{\alpha}_{t,\varepsilon}(y))^n\right|\\
&\leq c\,|x-y|^\lambda\int_{[0,t]^{2n}}\int_{\mathbb{R}^{n}}\prod_{j=1}^n|p_{j}|^{\lambda}
\exp\prod_{j_1=1}^\infty(1+4\lambda_{j_1}^2)^{-\frac14}\prod_{j_2=1}^\infty(1+4\widetilde{\lambda}_{j_2}^2)^{-\frac14}dpdrds\\
&\leq c_n\, \int_{\widetilde{D}^{2n}}\left(\int_{\mathbb{R}^n}\prod_{j=1}^n|p_{j}|^{\lambda}\prod_{j_1=1}^\infty(1+4\lambda_{j_1}^2)^{-\frac14}dp\right)^{\frac{H_2}{H_1+H_2}}\\
&\qquad\times\left(\int_{\mathbb{R}^n}\prod_{j=1}^n|p_{j}|^{\lambda}\prod_{j_2=1}^\infty(1+4\widetilde{\lambda}_{j_2}^2)^{-\frac14}dp\right)^{\frac{H_1}{H_1+H_2}}drds,
\end{align*}
where $\widetilde{D}=\{0<s_1<s_2<\cdots<s_n<t\}$.

Making the change of variables $\xi_j=\sum_{i=j}^np_i$ for $j=1,2,\cdots, n$ with the convention $\xi_{n+1}=0$, then by \eqref{sec3-DILT-eq.3} and \eqref{sec3-DILT-eq.ref},
\begin{align}\label{sec3-DILT-eq.7-1}
\int_{\mathbb{R}^n}\prod_{i=1}^n|p_{i}|^{\lambda}\prod_{j_1=1}^\infty(1+4\lambda_{j_1}^2)^{-\frac14}dp
&\leq c\,\int_{\mathbb{R}^n}\prod_{j=1}^n|\xi_{j}-\xi_{j+1}|^{\lambda}e^{-c_{H_1}\max_{1\leq j\leq n}|s_j-s_{j-1}||\xi_j|^{\frac1{H_1}}}d\xi\nonumber\\
&\leq c\,\sum_{S_2}\int_{\mathbb{R}^n}\prod_{j=1}^n|\xi_{j}|^{\lambda\alpha_j}e^{-c_{H_1}\max_{1\leq j\leq n}|s_j-s_{j-1}||\xi_j|^{\frac1{H_1}}}d\xi\nonumber\\
&\leq c\,\int_{\mathbb{R}^n}\prod_{j=1}^n|\xi_{j}|^{\lambda\alpha_j}\left(e^{-c_{H_1}\sum_{j=1}^n|s_j-s_{j-1}||\xi_j|^{\frac1{H_1}}}\right)^{\frac1n}d\xi\nonumber\\
&\leq c_{H_1,n}\prod_{j=1}^n|s_j-s_{j-1}|^{-H_1(\lambda\alpha_j+1)},
\end{align}
where in the second inequality we use
\begin{align*}
\prod_{j=1}^{n}|\xi_{j}-\xi_{j+1}|^{\lambda} &\leq c\,\prod_{j=1}^{n} (|\xi_{j}|^{k}+|\xi_{j+1}|^{\lambda})\\
&=\sum_{S_1}\prod_{j=1}^{n}(|\xi_{j}|^{\lambda\gamma_{j}}|\xi_{j+1}|^{\lambda\overline{\gamma_{j}}})\\
&\leq\sum_{S_2}\prod_{j=1}^{n}(|\xi_{j}|^{\lambda\alpha_{j}}),
\end{align*}
with $$S_1=\left\{\gamma_{j}, ~\overline{\gamma_{j}}: ~\gamma_{j}\in\{0,1\}, ~\gamma_{j}+\overline{\gamma_{j}}=1,  ~j=1, \cdots, n\right\}$$
and
$$
S_2=\left\{\alpha_{j}: ~\alpha_{j}\in\{0, 1, 2\}, ~j=1, \cdots n\right\}.
$$

Similarly, we have the same inequality for $\widetilde{\lambda}_{j_2}$ as \eqref{sec3-DILT-eq.7-1}. Thus,
\begin{align*}
\left|\mathbb{E}\left(\widetilde{\alpha}_{t,\varepsilon}\right)^n\right|
&\leq c_n|x-y|^{\lambda}\left(\sum_{S_2}\int_{\{0<r_1<\cdots<r_n<t\}}\prod_{j=1}^n(r_j-r_{j-1})^{-\frac{H_1H_2}{H_1+H_2}(1+\lambda\alpha_j)}dr\right)^2.
\end{align*}
So,
we need the condition $\frac{H_1H_2}{H_1+H_2}(2\lambda+1)<1$ to make sure
\begin{align}\label{sec3-DILT-eq.8}
&\left|\mathbb{E}(\widetilde{\alpha}_{t}(x)-\widetilde{\alpha}_{t}(y))^n\right|\nonumber\\
&\leq c_n\left|\mathbb{E}(\widetilde{\alpha}_{t}(x)-\widetilde{\alpha}_{t,\varepsilon}(x))^n\right|
+c_n\left|\mathbb{E}(\widetilde{\alpha}_{t,\varepsilon}(x)-\widetilde{\alpha}_{t,\varepsilon}(y))^n\right|\nonumber\\
&\qquad+c_n\left|\mathbb{E}(\widetilde{\alpha}_{t,\varepsilon}(y)-\widetilde{\alpha}_{t}(y))^n\right|\nonumber\\
&\leq c_{n,H_1,H_2,t}\,|x-y|^{n\lambda}.
\end{align}

Hence, from the Kolmogorov continuity criterion, the H\"{o}lder continuity of $\widetilde{\alpha}_{t}(y)$ in space variable of any order $\widetilde{\lambda}$ strictly less than $\min\{1,\frac{H_1+H_2}{2H_1H_2}-\frac{1}2\}$.

\textbf{H\"{o}lder continuity of time variable.}
Without loss of generality, we assume that $t<\widetilde{t}$. Then
\begin{align*}
&\left|\mathbb{E}(\widetilde{\alpha}_{t,\varepsilon}(y)-\widetilde{\alpha}_{\widetilde{t},\varepsilon}(y))^n\right|\\
&\leq\frac1{(2\pi)^{n}}\left(\int_{[t,\widetilde{t}]^{2n}}+\int_{[0,t]^n\times[t,\widetilde{t}]^n}\right)\int_{\mathbb{R}^{n}}
\prod_{j_1=1}^\infty(1+4\lambda_{j_1}^2)^{-\frac14}\prod_{j_2=1}^\infty(1+4\widetilde{\lambda}_{j_2}^2)^{-\frac14}dpdrds\\
&=:II_1+II_2.
\end{align*}

By \eqref{sec3-DILT-eq.7}, we have
\begin{align*}
II_1
&\leq c_n\left(\int_{\{t<r_1<\cdots<r_n<\widetilde{t}\}}\prod_{j=1}^n(r_j-r_{j-1})^{-\frac{H_1H_2}{H_1+H_2}}dr\right)^2\\
&\leq c_n\left((\widetilde{t}-t)^{\sum_{j=1}^n\left(1-\frac{H_1H_2}{H_1+H_2}\right)}\right)^2\\
&\leq c_n(\widetilde{t}-t)^{2n\left(1-\frac{H_1H_2}{H_1+H_2}\right)}
\end{align*}
and similarly,
\begin{align*}
II_2
\leq c_n(\widetilde{t}-t)^{n\left(1-\frac{H_1H_2}{H_1+H_2}\right)}.
\end{align*}
This gives
\begin{align}\label{sec3-DILT-eq.9}
&\left|\mathbb{E}(\widetilde{\alpha}_{t}(x)-\widetilde{\alpha}_{\widetilde{t}}(x))^n\right|\nonumber\\
&\leq c_n\left|\mathbb{E}(\widetilde{\alpha}_{t}(x)-\widetilde{\alpha}_{t,\varepsilon}(x))^n\right|
+c_n\left|\mathbb{E}(\widetilde{\alpha}_{t,\varepsilon}(x)-\widetilde{\alpha}_{\widetilde{t},\varepsilon}(x))^n\right|\nonumber\\
&\qquad+c_n\left|\mathbb{E}(\widetilde{\alpha}_{\widetilde{t},\varepsilon}(x)-\widetilde{\alpha}_{\widetilde{t}}(x))^n\right|\nonumber\\
&\leq c_{n,H_1,H_2,t}\,|\widetilde{t}-t|^{n\widetilde{\beta}}
\end{align}
with $\widetilde{\beta}<1-\frac{H_1H_2}{H_1+H_2}$.  This completes the proof.

\subsection{The proof of Theorem \ref{sec1-DCLT}}\label{sec3.3}

In this section we will proof the existence and H\"{o}lder continuity for the case of CLT. Compared with cases of SLT
and ILT, integral structure of CLT here is relatively simple. Then we can use the ways of SLT case directly here.

\textbf{Existence in $L^p$.}
By the definition of CLT in \eqref{sec1-eq.DCLT}, for any integer $n\geq1$,
\begin{align*}
\left|\mathbb{E}(\overline{\alpha}_{t,\varepsilon}(y))^n\right|
&\leq\frac1{(2\pi)^{n}}\int_{[0,t]^n}\int_{\mathbb{R}^{n}}e^{-\frac{\varepsilon}{2}\sum_{j=1}^n|p_j|^2}
\left|\mathbb{E}e^{\iota \sum_{i=1}^np_i(X^{H_1}_s-\widetilde{X}^{H_2}_s)}\right|dpds\\
&\leq\frac{1}{(2\pi)^n}\int_{[0,t]^{n}}\int_{\mathbb{R}^n}\prod_{j_1=1}^\infty(1+4\lambda_{j_1}^2)^{-\frac14}
\prod_{j_2=1}^\infty(1+4\widetilde{\lambda}_{j_2}^2)^{-\frac14}dpds.
\end{align*}

Using the same method as in \eqref{sec3-DILT-eq.5}, we have
\begin{align}\label{sec3-DCLT-eq.1}
&\int_{\mathbb{R}^n}\prod_{j_1=1}^\infty(1+4\lambda_{j_1}^2)^{-\frac14}\prod_{j_2=1}^\infty(1+4\widetilde{\lambda}_{j_2}^2)^{-\frac14}dp\nonumber\\
&\leq c\,\left(\int_{\mathbb{R}^n}\prod_{j_1=1}^\infty(1+4\lambda_{j_1}^2)^{-\frac14}dp\right)^{\frac{H_2}{H_1+H_2}}
\left(\int_{\mathbb{R}^n}\prod_{j_2=1}^\infty(1+4\widetilde{\lambda}_{j_2}^2)^{-\frac14}dp\right)^{\frac{H_1}{H_1+H_2}}\nonumber\\
&\leq c\,\prod_{j=1}^n|s_j-s_{j-1}|^{-\frac{2H_1H_2}{H_1+H_2}}.
\end{align}

Thus,
\begin{align*}
\left|\mathbb{E}\left(\overline{\alpha}_{t,\varepsilon}\right)^n\right|
&\leq c\,\int_{[0,t]^n}\prod_{j=1}^n|s_j-s_{j-1}|^{-\frac{2H_1H_2}{H_1+H_2}}ds\\
&\leq c\int_{\{0<s_1<\cdots<s_n<t\}}\prod_{j=1}^n|s_j-s_{j-1}|^{-\frac{2H_1H_2}{H_1+H_2}}ds,
\end{align*}
which is finite under condition $\frac{H_1H_2}{H_1+H_2}<\frac12$.
It is easy to see that $\left\{\overline{\alpha}_{t,\varepsilon}, \varepsilon>0\right\}$ is a Cauchy sequence and we get $\overline{\alpha}_{t,\varepsilon}$ exists in $L^p$, for all $p\geq1$, under the condition $\frac{H_1H_2}{H_1+H_2}<\frac12$.

\textbf{H\"{o}lder continuity of space and time variables.}

According to eq. \eqref{sec3-DCLT-eq.1}, we find that the result of H\"{o}lder continuity can be obtained by replacing $\frac{H_1H_2}{H_1+H_2}$ with $\frac{2H_1H_2}{H_1+H_2}$ in the case of ILT. Then similar to \eqref{sec3-DILT-eq.8} and \eqref{sec3-DILT-eq.9}, we can obtain
\begin{align*}
\left|\mathbb{E}(\overline{\alpha}_{t}(x)-\overline{\alpha}_{t}(y))^n\right|
\leq c_{n,H_1,H_2,t}\,|x-y|^{n\overline{\lambda}}
\end{align*}
for $\overline{\lambda}<\min\{1,\frac{H_1+H_2}{4H_1H_2}-\frac{1}2\}$; and
\begin{align*}
\left|\mathbb{E}(\overline{\alpha}_{t}(x)-\overline{\alpha}_{\widetilde{t}}(x))^n\right|\leq c_{n,H_1,H_2,t}\,|\widetilde{t}-t|^{n\overline{\beta}}
\end{align*}
with $\overline{\beta}<1-\frac{2H_1H_2}{H_1+H_2}$.  This completes the proof.

\bigskip

$\begin{array}{cc}
\begin{minipage}[t]{1\textwidth}
{\bf  Qian Yu }\\
School of Statistics, East China Normal University, Shanghai 200241, China \\
\texttt{qyumath@163.com}
\end{minipage}
\hfill
\end{array}$

\bigskip

$\begin{array}{cc}
\begin{minipage}[t]{1\textwidth}
{\bf Guangjun Shen}\\
Department of Mathematics, Anhui Normal University, Wuhu 241000, China\\
\texttt{gjshen@163.com}
\end{minipage}
\hfill
\end{array}$

\bigskip

$\begin{array}{cc}
\begin{minipage}[t]{1\textwidth}
{\bf Xiuwei Yin}\\
Department of Mathematics, Anhui Normal University, Wuhu 241000, China\\
\texttt{xweiyin@163.com}
\end{minipage}
\hfill
\end{array}$

\end{document}